\documentclass[12pt]{amsart}
\begin{document}

\def\1#1{\overline{#1}}
\def\2#1{\widetilde{#1}}
\def\3#1{\widehat{#1}}
\def\4#1{\mathbb{#1}}
\def\5#1{\frak{#1}}
\def\6#1{{\mathcal{#1}}}

\def\C{{\4C}}
\def\R{{\4R}}
\def\N{{\4N}}
\def\Z{{\4Z}}

\def\cp#1{\C{\rm P}^{#1}}

\def\frac#1#2{{#1\over#2}}
\def\square{\sqcup\!\!\!\!\sqcap}
\def\sec{'\! '}
\def\limind{\mathop{\oalign{lim\cr\hidewidth$\longrightarrow$\hidewidth\ 
cr}}}
\def\norm{\mid\!\mid}

\title[On Levi-flat hypersurfaces with prescribed boundary]{On Levi-flat hypersurfaces with prescribed boundary}
\author[P. Dolbeault, G. Tomassini and D. Zaitsev]{Pierre Dolbeault, Giuseppe Tomassini and Dmitri Zaitsev}
\address{Giuseppe Tomassini: Scuola Normale Superiore, Piazza dei Cavalieri 7, 56126 Pisa, ITALY}
\email{g.tomassini@sns.it}
\address{Pierre Dolbeault: Institut de Math\'ematiques de Jussieu, Universit\'e
Pierre et Marie Curie, 175 rue du Chevaleret, 75013 Paris, FRANCE}
\email{pierre.dolbeault@upmc.fr}
\address{Dmitri Zaitsev: School of Mathematics, Trinity College Dublin, Dublin 2, Ireland}
\email{zaitsev@maths.tcd.ie}
%\subjclass{}
\dedicatory {Dedicated to Professor J. J. Kohn on the occasion of his $75$th birthday}
\thanks{The research of the third author was supported in part by the Science Foundation Ireland grant 06/RFP/MAT018}
\maketitle
%\tableofcontents

\begin{abstract}
We address the problem of existence and uniqueness of a Levi-flat
hypersurface $M$ in $\C^n$ with prescribed compact boundary 
$S$ for $n\ge3$.
The situation for $n\ge3$ differs sharply from 
the well studied case $n=2$.
We first establish necessary conditions on $S$ at both complex and 
CR points, needed for the existence of $M$.
All CR points have to be nonminimal and all 
complex points have to be ``flat''.
Then, adding a positivity condition at complex points,
which is similar to the ellipticity for $n=2$
and excluding the possibility of $S$ to contain complex
$(n-2)$-dimensional submanifolds,
we obtain a solution $M$ to the above problem as
a projection of a possibly singular Levi-flat hypersurface
in $\R\times\C^n$. It turns out that $S$ has to be
a topological sphere with two complex points
and with compact CR orbits, also topological spheres, 
serving as boundaries of 
the (possibly singular) complex leaves of $M$.
There are no more global assumptions on $S$
like being contained in the boundary of
a strongly pseudoconvex domain, as it was in case $n=2$.
Furthermore, we show in our situation that any other
Levi-flat hypersurface with boundary $S$
must coincide with the constructed solution.
\end{abstract}

%\contentsline {section}{\tocsection {}{1}{Introduction}}{3}
%\contentsline {section}{\tocsection {}{2}{Local analysis and flatness conditions.}}{5}
%\contentsline {section}{\tocsection {}{3}{Global consequences of the local flatness}}{12}
%\contentsline {section}{\tocsection {}{4}{On boundaries of families of holomorphic chains
%with $C^\infty$ parameters.}}{15}
%\contentsline {section}{\tocsection {}{5}{On some Levi-flat $(2n-1)$-subvarieties with given
%boundary in ${\mathbb {C}}^n$.}}{26}
%\contentsline {section}{\tocsection {}{6}{About regularity.}}{29}
%\contentsline {section}{\tocsection {}{}{References}}{37}

%\def\Label#1{\label{#1}{\bf (#1)}~}
\def\Label#1{\label{#1}}

% Standard sets

\def\cn{{\C^n}}
\def\cnn{{\C^{n'}}}
\def\ocn{\2{\C^n}}
\def\ocnn{\2{\C^{n'}}}

% Abbreviations

\def\const{{\rm const}}
\def\rk{{\rm rank\,}}
\def\id{{\sf id}}
\def\aut{{\sf aut}}
\def\Aut{{\sf Aut}}
\def\CR{{\rm CR}}
\def\GL{{\sf GL}}
\def\Re{{\sf Re}\,}
\def\Im{{\sf Im}\,}

\def\codim{{\rm codim}}
\def\crd{\dim_{{\rm CR}}}
\def\crc{{\rm codim_{CR}}}

\def\phi{\varphi}
\def\eps{\varepsilon}
\def\d{\partial}
\def\a{\alpha}
\def\b{\beta}
\def\g{\gamma}
\def\G{\Gamma}
\def\D{\Delta}
\def\Om{\Omega}
\def\k{\kappa}
\def\l{\lambda}
\def\L{\Lambda}
\def\z{{\bar z}}
\def\w{{\bar w}}
\def\t{\tau}
\def\th{\theta}
\def\sec{''}
\def\cp#1{\C{\rm P}^{#1}}

\emergencystretch15pt
\frenchspacing

\newtheorem{Thm}{Theorem}[section]
\newtheorem{Cor}[Thm]{Corollary}
\newtheorem{Pro}[Thm]{Proposition}
\newtheorem{Lem}[Thm]{Lemma}
\newtheorem{Exa}[Thm]{Example}
\theoremstyle{definition}\newtheorem{Def}[Thm]{Definition}
\theoremstyle{remark}\newtheorem{Rem}[Thm]{\bf Remark}

\def\bpf{\begin{Dim}}
\def\bl{\begin{Lem}}
\def\el{\end{Lem}}
\def\bp{\begin{Pro}}
\def\ep{\end{Pro}}
\def\bt{\begin{Thm}}
\def\et{\end{Thm}}
\def\bc{\begin{Cor}}
\def\ec{\end{Cor}}
\def\bd{\begin{Def}}
\def\ed{\end{Def}}
\def\br{\begin{Rem}}
\def\er{\end{Rem}}
\def\be{\begin{Exa}}
\def\ee{\end{Exa}}
\def\bpf{\begin{proof}}
\def\epf{\end{proof}}
\def\ben{\begin{enumerate}}
\def\een{\end{enumerate}}

%\noindent\today

\section{Introduction}
Let $S$ be a smooth $2$-codimensional real submanifold of $\C^n$, $n\ge 2$.The problem of
finding a Levi-flat hypersurface $M\subset\C^n$  with
boundary $S$ has been extensively studied for $n=2$ by the methods of geometric theory of
several complex variables (cf. [BeG], [BeK], [Sh], [Kr], [ChS], [SlT] and  [ShT] for the
unbounded case). The starting point of Bedford and Gaveau [BeG] is the classical theorem of
Bishop on the existence of local 1-parameter families of analytic discs [Bi].
Near an elliptic complex point $p\in S$,
the real surface $S\setminus \{p\}\subset\C^2$ is foliated by smooth compact real curves which bound
analytic discs in $\C^2$. The family of these discs fills a smooth
Levi-flat hypersurface. Assume that $S$ is contained in a strictly pseudoconvex compact
hypersurface, and it is the graph of a smooth function $g$ on the sphere
$S^2\subset\C\times\R$ with only two complex points that are both elliptic. Then in [BeG] 
it is proved that the families of analytic discs near complex
points  extend to one global family filling a topological 3-dimensional ball $M$ bounded by
$S$.  Moreover, $M$ is the hull of holomorphy of $S$. The more general case, when $S$ is
contained in the boundary of a bounded strictly pseudoconvex domain of a $2$-dimensional
Stein manifold have been studied in [BeK] and [Kr].

In this paper we address the corresponding problem of finding Levi-flat hypersurfaces with prescribed compact boundary in $\C^n$ for $n>2$. The
situation here turns out to be totally different from what it is in $\C^2$. The first difference is already visible from the elementary observation that a submanifold of real codimension $2$ in general position is totally real in $\C^2$ but no more such in $\C^n$ with $n>2$.
Furthermore, if a surface $S\subset\C^2$ is real-analytic, any real-analytic foliation of
$S$ by real curves extends locally to a foliation of $\C^2$ by complex curves and hence $S$ locally bounds many possible Levi-flat hypersurfaces.
On the other hand, in higher dimension, a real-analytic 
submanifold $S\subset\C^n$ of codimension $2$ in general position
does not even locally bound a Levi-flat hypersurface $M$. Indeed,
locally $S$ is the graph of
a smooth function $g$ over a real hypersurface in $\C^{n-1}$, 
so the existence of a local Levi-flat graph
extending $S$ amounts to solving a boundary problem for a system of
three (quasi-linear, elliptic degenerate) operators and this requires
nontrivial compatibility conditions for $g$.
In other words, the problem we are dealing with becomes overdetermined. 

In Section 2 we study the necessary 
local compatibility conditions needed
for a $2$-codimensional smooth real submanifold $S\subset\C^n$ 
to bound a Levi-flat hypersurface at least locally.
First we observe that near a CR point, $S$ must be
nowhere minimal (in the sense of Tumanov [Tu], 
see \S2.1), i.e. all local CR orbits must be of
positive codimension. Next we consider a complex point $p\in S$ and local holomorphic
coordinates  
$(z,w)\in \C^{n-1}\times \C$, vanishing at $p$, such
that $S$ is locally given by the equation
\begin{equation} 
w= Q(z)+O(|z|^3),
\end{equation}
where $ Q(z)$ is a complex valued quadratic form in the real coordinates
$(\Re z, \Im z)\in\R^{n-1}\times\R^{n-1}$ or, equivalently, in $(z,{\bar{z}})$.
We observe that not all quadratic forms $Q$ can appear 
when $S$ bounds a Levi-flat hypersurface.
In fact, $Q$ has to satisfy a certain flatness condition,
in which case we call $p$ ``flat''.
We further call a flat point $p\in S$ ``elliptic'' if $Q(z)\in {\mathbb R}_+$  for every $z\ne0$ in suitable coordinates (cf. Definition~\ref{elliptic}). Elliptic flat complex points are always isolated
(see Remark~\ref{isolated}).

In our main result, we take these local necessary conditions 
as our assumptions and obtain a global conclusion.
There are no more global assumptions on $S$ (like being
contained in the boundary of a strongly pseudoconvex domain) 
as it was the case in $\C^2$.
Each of our assumptions can be checked in a neighborhood of a point and is therefore purely local.
There is a technical difficulty, however, that the possible Levi-flat hypersurface
solving the boundary problem may have self-intersections
producing singularities even along large sets.
We illustrate this phenomenon by the following example:

\be
Let $\4S^{2n-2}\subset \C^{n-1}\times\R\cong \R^{2n-1}$ be the unit sphere
and consider an immersion $\g\colon [-1,1]\to \C$ with $\g(x)\ne \g(-x)$ for all $x\ne 0$.
Then it is easy to see that the image $S\subset \C^{n}$ of $\4S^{2n-2}$
under the map $\id\times\g\colon \C^{n-1}\times\R\to \C^{n-1}\times\C$
is a submanifold of codimension $2$ that locally bounds an
immersed Levi-flat hypersurface
obtained as the image under the same map of the unit ball $\4B^{2n-1}$ 
bounded by $\4S^{2n-2}$ in $\C^n$. However, the image of the ball itself
can be singular at points $(z,\g(x))$ if there exist points $x'\ne x$ with $\g(x)=\g(x')$.
Clearly the set of such points can be very large.
\ee

In view of this example we have to allow more general Levi-flat ``hypersurfaces''
that are obtained as images of real manifolds.
Furthermore, the complex leaves themselves may have
singularities away from their boundaries.
Hence we are led to the following real analogue 
of complex-analytic varieties (cf.\ [HaL1]):

\bd\Label{subvar}
A closed subset $M$ in a domain $\Omega\subset\C^n$ is said to be a
$d$-{\it subvariety (with negligible singularities)},
if it is of locally finite  Hausdorff $d$-dimensional measure
and there exists a closed subset $\sigma$ of $M$ of $d$-dimensional Hausdorff measure zero
such that $M\setminus\sigma$ is an oriented smooth real $d$-dimensional submanifold.
The minimal set $\sigma$ with this property is called the 
{\it singular set} ${\rm Sing}\,M$ and its complement in $M$ 
the {\it regular set} ${\rm Reg}\,M$.
$M$ is said to be {\it Levi-flat} if ${\rm Reg}\,M$ is a Levi-flat 
CR submanifold in the usual sense (cf.\ Section 2).
\ed

%that we call ``$(2n-1)$-subvarieties in the abstract sense'',
%where $(2n-1)$ indicates the dimension.
Finally, we add a local assumption
on $S$ guaranteeing that all the CR orbits have the same dimension
and hence define a smooth foliation away from the complex points. We have:

\bt\Label{main}
Let $S\subset \C^n$, $n>2$, be a compact connected smooth real $2$-codimensional submanifold
such that the following holds: 
\begin{enumerate}
\item[(i)] $S$ is nonminimal at every CR point;
\item[(ii)] every complex point of $S$ is flat and elliptic and there exists at least one
such point; 
\item[(iii)] $S$ does not contain complex submanifolds of dimension $n-2$. 
\end{enumerate}
Then $S$ is a topological sphere with two complex points and 
there exist a smooth submanifold $\2S$
and a Levi-flat $(2n-1)$-subvariety $\2M$ in $\R\times\C^n$ (i.e. $\2M$ is Levi-flat in $\C\times\C^n$),
both contained in $[0,1]\times\C^n$
such that $\2S=d\2M$ in the sense of currents
and the natural projection $\pi\colon [0,1]\times\C^n\to \C^n$
restricts to a diffeomorphism between $\2S$ and $S$.
\et

As mentioned before, (i) and the flatness in (ii) are necessary
for $S$ to locally bound a Levi-flat hypersurface.
The condition (iii) may appear artificial at a first glance.
However, J. Lebl \cite{L} has constructed an example
of a generic nowhere minimal real-analytic submanifold 
of codimension $2$ in $\C^3$ that does not satisfy (iii) and
that locally does not bound any smooth Levi-flat hypersurface.

In addition, we obtain the following more precise information:

\bt\Label{precise}
The Levi-flat $(2n-1)$-subvariety $\2M$ in Theorem~\ref{main} can be
chosen with the following properties:
\begin{enumerate}
\item[(i)] $\2S$ has two complex points $p_0$ and $p_1$
with $\2S\cap(\{j\}\times\C^n)=\{p_j\}$ for $j=0,1$;  every other slice $\{x\}\times\C^n$ with $x\in(0,1)$, intersects $\2S$ transversally along a submanifold diffeomorphic to a sphere that bounds (in the sense of currents) the (possibly singular) irreducible complex-analytic hypersurface $(\2M\setminus\2S)\cap (\{x\}\times\C^n)$;
\item[(ii)] the singular set ${\rm Sing}\,\2M$ is the union of $\2S$ and a 
closed subset of $\2M\setminus\2S$ 
of Hausdorff dimension at most $2n-3$; moreover each slice
$({\rm Sing}\,\2M\setminus\2S)\cap (\{x\}\times\C^n)$
is of Hausdorff dimension at most $2n-4$;
\item[(iii)] there exists a closed subset $\2A\subset\2S$ 
of Hausdorff $(2n-2)$-dimensional measure zero
such that away from $\2A$, $\2M$ is a smooth submanifold with 
boundary $\2S$ near $\2S$;
moreover $\2A$ can be chosen such that each slice $\2A\cap (\{x\}\times\C^n)$
is of Hausdorff $(2n-3)$-dimensional measure zero.
\end{enumerate}
\et

Finally, we obtain the following uniqueness result:

\bt\Label{unique}
Let $S$, $\2S$ and $\2M$ satisfy the conclusions of Theorems~\ref{main}
and \ref{precise}.
Suppose that $M\subset\C^n$ 
is another Levi-flat $(2n-1)$-subvariety
such that the following holds:
\begin{enumerate}
\item[(i)] the singular set of $M$ is the union of $S$ and 
a closed subset of $M\setminus S$ of Hausdorff dimension at most $2n-3$;
\item[(ii)] through every point in $M\setminus S$ there is a complex-analytic subset of $M$ of positive dimension;
\item[(iii)] there exists a closed subset $A\subset S$ of Hausdorff dimension at most $2n-3$ such that, away from $A$, $M$ is a manifold with boundary $S$ near $S$, and such that
for every leaf $L\subset {\rm Reg}\,M$, 
one has $\1L\cap S\not\subset A$.
\end{enumerate}
Then $M$ coincides with the image $\pi(\2M)$.
\et

Here by a leaf in the Levi-flat hypersurface ${\rm Reg}\,M$ 
we mean a maximal image of an injective holomorphic immersion 
of a connected complex $(n-1)$-dimensional manifold.

The first main step of the proof, given in Section 2 and 3,
consists of proving the following:

\bp
Under the assumptions of Theorem~\ref{main}, 
$S$ is diffeomorphic to the unit
sphere $\4S^{2n-2}\subset \C^{n-1}_z\times\R_x$; 
there are precisely two complex points of $S$ corresponding to the poles $\{x=\pm 1\}$;
the CR orbits of $S$ away from the poles correspond to the spheres
$\{x={\rm const}\}$. 
\ep

This result allows us to solve the boundary problem mentioned above by applying
a theorem of Harvey-Lawson [HaL1] leaf-wise for single CR orbits.
Thus, we are led to prove the version of the latter theorem with smooth
parameters (cf.\ Theorem 4.1) which also extends results of [Do1] 
from the real-analytic to the smooth case (cf. [Do1, Th\'eor\`eme 6.9]). 
Next we prove that the CR orbits can be represented as the level surfaces of a smooth function $\nu$
having a smooth extension to the complex points.
Consequently, the graph $N$ of $\nu$ is a submanifold of the real
hyperplane
$E := \R_x\times \C^{n}$ of $\C\times\C^{n}$.
Moreover the CR orbits are sent by $\nu\times \id$ into the
hyperplane sections $N\cap \{x={\rm  const}\}$.  
Then we shall derive from Theorem 4.1 that 
there exists a unique Levi-flat $(2n-1)$-subvariety $\2M\subset E$, 
foliated by complex $(n-1)$-subvarieties, such that
$d\2M=N$ in the sense of currents, 
leading to the conclusion of Theorem~\ref{main}.

In conclusion, we mention another related problem, suggested by Gaveau (cf. [Ga]),
of finding the Levi-degenerate hypersurfaces with prescribed boundary 
(instead of the Levi-flat ones). In this context the boundary problem was studied in [SlT].

Some results of this paper have been announced in \cite{DTZ}.

\vskip 2mm
{\bf Acknowledgments.} From P.~D. to G. Henkin for conversations at the 
beginning of this work. From P.~D. and G.~T. to the European networks "Analyse Complexe et 
G\'eom\'etrie Analytique" for support.
From D.~Z. for hospitality of the Scuola Normale Superiore during his stay in Pisa.

\section{Local analysis and flatness conditions.}
We recall some definitions and fix the notation. 
Let $M$ be a smooth, connected, real
submanifold of $\C^n$, $n\ge 2$, of real dimension $m$. For every $p\in M$ we denote
$H_p(M)=T_p(M)\cap iT_p(M)$ the complex tangent space of $M$ at $p$ and $HM$ the fiber
bundle $\cup_{p\in M} H_p(M)$. ${\rm dim}_{_{CR}}\,M_p:={\rm dim}_{_{\C}}\,H_pM$ is
called the CR {\it dimension} of $M$ at $p$, ${\rm codim}_{_{CR}}\,M_p:=m-2 {\rm
dim}_{_{\C}}\,H_pM$ the CR {\it codimension}. $M$ is said to be a $CR$ {\it manifold} if
${\rm dim}_{_{\C}}\,H_p M$ is constant. In this case  ${\rm dim}_{_{CR}} M:={\rm
dim}_{_{\C}}\,H_pM$ is, by definition, the CR {\it dimension} of $M$, and 
${\rm codim}_{_{CR}}\,M:=m-2{\rm dim}_{_{CR}}
M$ the CR {\it codimension}. The CR dimension satisfies $m-n\le {\rm dim}_{_{CR}}\,M\le
m/2$. A CR  submanifold $M\subset \C^n$ of CR dimension $m-n$ (i.e. such that $T_pM+iT_pM=\C^n$)
is called {\it generic}. A CR submanifold $M\subset \C^n$ is called {\it maximally complex}
if $m=2{\rm dim}_{_{CR}}M+1$. Finally we recall that a $CR$ submanifold $M$ is said to be
{\it minimal} (in the sense of Tumanov [Tu]) at a point $p$ if there does not exist a submanifold
$N$ of $M$ of lower dimension through $p$ such that $HN = HM|_{N}$. In the general case, 
by a theorem  of Sussman (and a theorem of Nagano in the real-analytic case),
all possible submanifolds $N$ with this property contain, as germs at $p$, one of the
minimal possible dimension, called a (local) CR {\it orbit} of $p$ in $M$ (cf.
e.g. \cite{BERbook}). The germ at $p$ of the CR orbit of $p$ is uniquely
determined. 
In particular, $M$ is minimal at
$p$ if and only if a CR orbit of $p$ is open and hence, of the maximal possible dimension.
On the other hand, the minimal possible dimension of a CR orbit is $\dim_\R H_pM=2l$
which is twice the CR dimension of $M$. It is easy to see that all CR orbits in $M$ have
their dimension equal to $2l$  if and only if $M$ is Levi-flat i.e. the Levi form of $M$
identically vanishes.
\subsection{Behaviour near CR points.}
Let $S$ be a smooth real submanifold of real codimension $2$ in $\C^n$ (not  
necessarily compact). We say that $S$ is a {\it locally flat boundary} at a point $p\in S$
if there exist an open neighbourhood $U$ of $p$ in $\C^n$ and a Levi-flat hypersurface
$M\subset U$ with boundary $U\cap S$. Assume that $S$ is a   locally flat boundary and let
$p\in S$ be such that $S$ is $CR$ near $p$. Then, near $p$, 
$S$ is either a complex hypersurface (in which case it is clearly a locally flat  
boundary) or a generic submanifold of $\C^n$ at least at some points. In the second case
being a locally flat boundary turns out to be a nontrivial condition for $n\ge 3$. Indeed, 
suppose that $M\subset\C^n$ is a Levi-flat hypersurface bounded by a generic submanifold
$S$. Consider the foliation by complex hypersurfaces of $M$ where it extends smoothly to the
boundary. Since the boundary $S$ is generic, it cannot be tangent to a complex leaf.
Hence the leaf $M_p$ of $M$ through $p$ intersects $S$ transversally along a real
hypersurface $S_p\subset S$.  Since $H_qS\subset H_qM = T_qM_p$ for $q\in S_p$ near $p$, it
follows that $H_qS_p=H_qS$ for such $q$. Hence $S$ cannot be minimal at $p$ with $p$ being
arbitrary generic smooth boundary point. In fact, it follows that $S_p$ is either a single 
$CR$ orbit of $S$ or a union of $CR$ orbits.

In case $n=2$, $S$ is totally real (i.e. $HS=\{0\}$)
in its CR points and hence is obviously nowhere
minimal. However, for $n\ge 3$, $S$ cannot be totally real for  reasons of dimension and its ``nowhere minimality'' becomes a nontrivial condition.
(Recall the standard fact that a general CR submanifold 
can always be locally perturbed to a Levi-nondegenerate one
and hence to a minimal one but not vice versa).

If all CR orbits of $S$ are $1$-codimensional, 
they can be locally graphed over minimal real hypersurfaces in $\C^{n-1}$.
As a consequence of the extension theorems of Tr\'epreau and Tumanov \cite{Tr86,Tu88}
(more precisely, of their parametric versions)
we obtain that the condition of nowhere minimality
is sufficient for $S$ to be a locally flat boundary.
On the other hand, if all CR orbits of $S$ are $2$-codimensional
and hence are complex submanifolds, then $S$ is
Levi-flat and is clearly a locally flat boundary.
In view of the condition (iii) in Theorem~\ref{main},
the latter situation does not occur in our case.

\subsection{Complex points and their fundamental forms.}
We now study the behaviour of a $2$-codimensional submanifold $S\subset\C^n$
near a complex point $p\in S$, i.e.\ $p$ is such that $T_pS$ is a complex hyperplane in $T_p\C^n$. 
In suitable holomorphic coordinates $(z,w)\in \C^{n-1}\times\C$ vanishing at $p$,
$S$ is locally given by an equation
\begin{equation}\Label{1}
w= Q(z) + O(|z|^3),\quad
Q(z)= \sum_{1\le i,j\le n-1} (a_{ij}z_iz_j + b_{ij}z_i\1z_j +c_{ij}\1z_i\1z_j),
\end{equation}
where $(a_{ij})$ and $(c_{ij})$ are symmetric complex matrices
and $(b_{ij})$ is an arbitrary complex matrix.
The form $Q(z)$ can be seen as a ``fundamental form'' of $S$ at $p$,
however it is not uniquely determined (even as a tensor).
In fact, a holomorphic quadratic change of coordinates 
of the form $(z,w)\mapsto (z,w+\sum a'_{ij}z_iz_j)$
results in adding $(a'_{ij})$ to the matrix $(a_{ij})$.
On the other hand, it can be easily seen that 
the matrices $(b_{ij})$ and $(c_{ij})$
transform as tensors $T_pS\times T_pS\to T_p\C^n/T_pS$
under all biholomorphic changes of $(z,w)$
preserving the form in \eqref{1}.
\subsection{A necessary condition at complex points.}
Clearly any symmetric $\C$-valued $\R$-bilinear form $Q$ on $\C^2$
can appear in the equation \eqref{1}. 
However we shall see that the condition
on $S$ to be a locally flat boundary at $p$ 
implies a nontrivial condition on $Q$.
\bd\Label{def-flat}
We call $S$ {\it flat} at a complex point $p\in S$ 
if, in some (and hence in any) coordinates $(z,w)$
as in \eqref{1}, the term $\sum b_{ij}z_i\1z_j$  
takes values in some real line in $\C$,
i.e. there exists a complex number $\l\in\C$ 
such that $\sum b_{ij}z_i\1z_j\in\l\R$ for all
$z=(z_1,\ldots,z_{n-1})$.
\ed

\be
If $(b_{ij})$ is a Hermitian matrix, $S$ as above is automatically flat at $p$.
Also in case $n=2$ the flatness always holds.
On the other hand, any submanifold $S\subset\C^3_{z_1,z_2,w}$ given by 
$w=|z_1|^2+i|z_2|^2+O(|z|^3)$ is not flat at $0$.
\ee
We have the following property:
\bl\Label{flat}
Let $S\subset\C^n$ be a locally flat boundary
with complex point $p\in S$.
Then $S$ is flat at $p$.
\el
\bpf
Suppose first that $S$ is contained in a Levi-flat hypersurface $M$.
If $M$ is real-analytic, we can choose local holomorphic coordinates $(z,w)$
such that $M$ is the hyperplane $\Im w=0$.
Since $S\subset M$, it follows from \eqref{1} that $Q(z)\in\R$ for all $z$.
Then $Q(z)+Q(iz)\in\R$ which implies, in view of
Definition~\ref{def-flat}, that $S$ is flat at $p$.
If $M$ is merely smooth,
it can be approximated by a real-analytic hyperplane up to order $3$ at $p$ which is
enough to conclude the flatness at $p$,
 as claimed.
Finally, the general case, when $S$ bounds a Levi-flat hypersurface $M$,
is reduced to the previous ones because $S$
can be approximated in the $C^2$ topology by a submanifold $\2S\subset M$
with a complex point.
\epf
% From Lemma 2.3 we immediately obtain:
%\bc
%If $S$ is a locally flat boundary near $p$,
%its fundamental class at $p$ is flat.
%\ec
%\bpf
%The statement follows by pushing the boundary $S$ inside the Levi-flat  
%piece $M$.
%\epf
If $S$ is flat, by making a change of coordinates $(z,w)\mapsto(z,\l w)$,
it is easy to make $\sum b_{ij}z_i\1z_j\in\R$ for all $z$.
Furthermore, by a change of coordinates $(z,w)\mapsto (z,w+\sum a'_{ij}z_iz_j)$ we can
choose the holomorphic term in \eqref{1} to be the conjugate of the antiholomorphic one and
so make the whole form $Q$ real-valued. Hence we can transform $Q$ to the following
normal form:
\bd\Label{normal}
We say that $S$ is in a {\em flat normal form} at $p$
if the coordinates $(z,w)$ as in \eqref{1} are chosen
such that $Q(z)\in\R$ for all $z\in\C^{n-1}$.
\ed
If the matrix $(b_{ij})$ is nonzero,
it follows from the discussion above
that in a flat normal form, 
$Q(z)$ is uniquely determined 
as a tensor $T_pS\times T_pS\to T_p\C^n/T_pS$.

\subsection{Elliptic flat points.}
We now assume the above necessary conditions for $S$ in order to be a locally flat boundary.
That is, we assume that $S$ is nonminimal in its generic points and flat at its
complex points. We then prove that $S$ is a locally flat boundary near $p$ assuming in
addition a certain positivity (ellipticity) condition. The latter condition is analogous to
that for $2$-surfaces in $\C^2$.
\bd\Label{elliptic}
Let $p\in S$ be a flat point. We say that $p$ is {\it elliptic} if, in some (and hence in
any) flat normal form, the real quadratic form $Q(z)$ is positive or negative definite.
\ed
By adding $Q(z)$ and $Q(iz)$
we see that  ellipticity implies that the matrix 
$(b_{ij})$ is positive definite.
\br
Definition~\ref{elliptic} generalizes the well-known notion of ellipticity
(in the sense of Bishop) for $n=2$. Note that $S\subset\C^2$ is always flat 
(in the sense of Definition~\ref{def-flat}) at any complex point. In general,
it can be shown that $S$ is elliptic at a flat complex point $p$
if and only if every intersection $L\cap S$
with a complex $2$-plane $L$ through $p$ satisfying $L\not\subset T_pS$
is elliptic in $L\cong\C^2$ in the sense of Bishop.
\er
\subsection{Quadrics with elliptic flat points}
The simplest example of $S$ is the quadric of $\C^3$:
\begin{equation}\Label{quad}
w=Q(z),
\end{equation}
where $Q$ is as in \eqref{1}.
In our case when $S$ is flat and elliptic at $p=0$,
we can choose the coordinates $(z,w)$
where $Q(z)$ is real and positive definite.
We have the following elementary properties
whose proof is left to the reader:

\bl\Label{2.7}
Suppose that the quadric \eqref{quad} is flat and elliptic at $0$.
Then it is CR and nowhere minimal outside $0$, and
the CR orbits are precisely the $3$-dimensional
ellipsoids given by $w=const$. The Levi form at the CR points
is positive definite.
\el

\br\Label{isolated}
In particular, it follows that 
elliptic flat points are always isolated complex points. 
This property also holds for general $2$-codimensional 
submanifolds $S\subset\cn$ as can be seen by comparing
$S$ with a corresponding approximating quadric.
\er

\subsection{Nowhere minimal surfaces with elliptic flat points}
For elliptic flat points, we now show that the above necessary conditions
are in some sense sufficient for $S$ to be a locally flat boundary.

\bp\Label{elliptic-solution}
Assume that $S\subset\cn$ ($n\ge 3$) is nowhere minimal at all its CR points
and has an elliptic flat complex point $p$. Then there exists a neighbourhood $V$ of $p$ such that $V\setminus \{p\}$ is foliated by compact real $(2n-3)$-dimensional CR orbits and there exists a smooth function $\nu$ without critical points away from $p$, having the CR orbits as its level surfaces.
\ep

\bpf
We may assume that $S$ is given by $w=\phi(z)$ with $p=0$ and $\phi(z)=Q(z)+O(|z|^3)$ as in \eqref{1}.
Denote by $S_0$ the corresponding quadric \eqref{quad} for $n\geq 3$.
%We write $\phi(z,\1z)$ (resp. $\phi_0(z,\1z)$) for the right-hand side
%of \eqref{1} (resp.\ of \eqref{quad}).
By differentiating \eqref{1} we obtain for the tangent spaces the 
following asymptotics:
\begin{equation}\Label{full-tang}
T_{(z,\phi(z))}S=T_{(z,Q(z))}S_0 + O(|z|^2), \quad z\in \C^{n-1},
\end{equation}
where the tangent spaces are understood as elements
of the Grassmannian of all $2$-codimensional real vector subspaces of $\R^{2n}$.
Here both $T_{(z,\phi(z))}S$ and $T_{(z,Q(z))}S_0$ 
depend continuously on $z$ near the origin.
Hence the choice of coordinates in the corresponding Grassmannian plays no role 
in \eqref{full-tang}. However, the corresponding complex tangent spaces
change their dimension at the origin and have there no limit in the Grassmannian
of all $2$-codimensional complex vector subspaces of $\C^{n}$.
Hence, some care is needed to choose suitable coordinates.
To do this, consider the unit ellipsoid $G:=\{z\in \C^{n-1} : Q(z)=1\}$ and the projection 
$$\pi\colon \C^{n}\setminus\{z=0\}\to G, \quad (z,w)\mapsto z/\sqrt{Q(z)}.$$
Then, for every $z\in G$, consider a real orthonormal basis $e_1(z),\ldots,e_{2n}(z)$
(with respect to the standard Euclidean product)
such that 
\begin{equation}\Label{basis-cond}
e_1(z),\ldots,e_{2n-4}(z)\in H_z G, \quad e_{2n-3}(z)\in T_z G.
\end{equation}
Locally such a basis can be chosen continuously depending on $z$.
For every $(z,w)$, $z\ne 0$, consider coordinate charts in the Grassmannians
of all vector subspaces of $T_{z,w}\C^n$ given by the basis
$e_1(\pi(z,w)),\ldots,e_{2n}(\pi(z,w))$. Then every vector subspace in such a chart is
represented by a matrix of suitable size. In particular, in view of \eqref{basis-cond},
the spaces $T_{(z,Q(z))}S_0$ and $H_{(z,Q(z))}S_0$ for $z\ne 0$ are represented by
the corresponding zero matrices.
A direct calculation using \eqref{full-tang} shows that, in the chosen coordinates,
\begin{equation}\Label{complex-tang}
H_{(z,\phi(z))}S=H_{(z,Q(z))}S_0 + O(|z|^2), \quad z\ne 0, \, z\to 0.
\end{equation}
In particular, $S$ is CR outside the origin.
%and the complex tangent space to $S$
%has the same asymptotics as in (5).
The next calculation of the Lie brackets of unit vector fields
in $HS_0$ and $HS$ shows that the Levi forms $L(S_0)$ and $L(S)$ are related by
\begin{equation}\Label{levi-rel}
L(S)_{(z,\phi(z))}=L(S_0)_{(z,Q(z))} + O(|z|), \quad z\ne 0, \, z\to 0.
\end{equation}
Taking into account that $L(S_0)_{(\l z,Q(\l z))}=\frac{1}{\l}L(S_0)_{(z,Q(z))}$
and normalizing the Levi form appropriately,
we obtain, for the normalized Levi forms,
\begin{equation}\Label{norm-levi-rel}
\frac{L(S)_{(z,\phi(z))}}{|L(S)_{(z,\phi(z))}|}
=\frac{L(S_0)_{(z,Q(z))}}{|L(S_0)_{(z,Q(z))}|} + O(|z|^2), \quad z\ne 0, \, z\to 0.
\end{equation}
In particular the Levi-form of $S$ is nonzero at CR points near $0$.

Denote by $E(q)$, $q\in S\setminus \{0\}$,
the tangent spaces to the local CR orbit of $S$ through $q$.
Since $S$ is nowhere minimal and the Levi form
is everywhere nonzero at CR points, the space $E(q)$, of
real dimension $2n-3$, is spanned by the complex tangent space $H_qS$ and
the Levi form at $q$. Denoting by $E_0(q)$, $q\in S_0$, the analogous object for  
$S_0$ and using the above asymptotics we obtain
\begin{equation}\Label{6}
E((z,\phi(z)))=E_0((z,Q(z))) + O(|z|^2), \quad z\ne 0, \quad z\to 0.
\end{equation}

We now show that CR orbits of $S$
have a transversality property with respect to the radial lines.
By Lemma \ref{2.7},
the CR orbits of $S_0$ are the ellipsoids
given by $w=\const$, i.e.\ in view of \eqref{quad}, by $Q(z)=c$, $w=c$, for  
$c\in\R_+$.
In particular, $d\pi$ projects each $E_0(q)$, $q\in S_0\setminus\{0\}$,
bijectively onto $T_{\pi(q)} G$.
We now conclude from \eqref{6} that the same property holds near the origin also
when $S_0$ is replaced by $S$ and $E_0$ by $E$. This is the crucial observation.
It implies that the restriction of $\pi$ to each (local) CR orbit of $S$
(in a suitable neighbourhood of the origin) is a local diffeomorphism.

Using the above observation we can now define global CR orbits.
We start with a point $q\in S\setminus\{0\}$
and extend the local inverse of the restriction of $\pi$ to 
the local CR orbit of $q$ in $S$ along paths on $G$.

We claim that
the extension is always possible if the length of the path is bounded
by a fixed constant and the starting point $q$ is sufficiently close to $0$
(depending on the constant).
Indeed, given the spaces $E_0$ and $E$, the CR orbits are obtained
by solving ODEs along the paths. 
We can regard $z=(z_1,\ldots,z_{n-1})$ as coordinates on both $S_0$ and $S$
and thus identify them locally with $\C^n$ by projecting to $z$-component.
In case of $S_0$, both CR orbits and ODEs are invariant under scaling 
$z\mapsto \lambda z$, $\lambda\in\R$. 
On the other hand, \eqref{6} implies that
the scaled ODEs for $S$ become arbitrarily close to those for $S$
provided $|z|$ is sufficiently small.
Since CR orbits of $S_0$ are the ellipsoids $Q(z)=\const$,
the claim follows from continuous dependence of
ODE solutions under perturbations.

Clearly there exists a constant $C>0$ 
such that we can reach any point of $G$ by a path of a length not exceeding $C$
and such that any two such paths are homotopic within
paths of length not exceeding $C>0$ (note that the ellipsoid $G$ is simply-connected for
$n\ge 3$). It follows that any two extensions of a CR orbit
along two paths as above coincide.
Hence, for every $q\in S$ sufficiently close to $0$, 
we obtain a global compact CR orbit $S_q$ through $q$
that projects diffeomorphically onto $G$ via $\pi$.

In order to obtain a function $\nu$ as desired,
consider a smooth curve $\gamma\colon \lbrack 0,\eps)\to S$
with $\g(0)=p$ that is diffeomorphic onto its image.
Define $\nu=\g^{-1}$ on the image of $\g$.
Then it follows from the above description that every CR orbit sufficiently close to $p$
intersects $\g(\lbrack 0,\eps))$ at precisely one point.
Hence there is a unique extension of $\nu$ from $\g(\lbrack 0,\eps))$
to a neighbourhood of $p$ having CR orbits as its level surfaces.
Obviously, $\nu$ is smooth away from $p$.
Furthermore, differentiating \eqref{6} and taking into account
the properties of $E_0((z,Q(z)))$, we conclude 
that directional derivatives of $E((z,\phi(z)))$ along the unit vectors
are $O(1/|z|)$. Since the diameter of the orbit $S_q$
is $O(|z|)$, we obtain that the derivative of $\nu$ near the origin is bounded.
Hence $\nu$ is Lipschitz up to the origin.
Furthermore, it follows from the construction that
all higher order derivatives of $\nu$ are $O(1/|z|^m)$ for suitable $m$ (depending on the derivative). Therefore
replacing $\nu$ by $e^{-1/\nu}$,
we may assume that $\nu$ is smooth also across $p$.
\epf

\br\Label{filling}
In the setting of Proposition~\ref{elliptic-solution}
one can construct a complex hypersurface $M_q\subset\C^n$
whose boundary is a global CR orbit $S_q$ in $S$ for $q\in S$
sufficiently close to the origin.
Indeed $S_q$ can be regarded
as the graph of a CR function $f$ over its projection
$\2S_q\subset \C^{n-1}_z$.
Then the graph of the unique holomorphic extension of $f$
to the interior of $\2S_q$ defines such a hypersurface $M_q$.
Both $S_q$ and $M_q$ are trivially projections 
of $\2S_q:=\{\nu(q)\}\times S_q$ and $\2M_q:=\{\nu(q)\}\times M_q$
under $\pi\colon \R\times\C^n\to\C^n$.
If $S_q$ is varying over all CR orbits near the origin,
the family of complex hypersurfaces $\2M_q$
fills a Levi-flat hypersurface $\2M\subset\R\times\C^n$
with boundary $\2S$ being the union of all $\2S_q$
and we obtain the conclusion of Theorem~\ref{main} near 
the given complex point $p$.
%Since $\2S_q$ is transversal to the
%(real) radial rays in $\C^{n-1}_z$,
%it follows that $\2S_{q_1}$ and $\2S_{q_2}$ are disjoint
%for $q_1,q_2\in S\setminus\{0\}$ from different CR orbits
%and the union of the hypersurfaces $M_q$ is 
%Then it is clear that there exists a continuous map $s$ as required.
\er

\section{Global consequences of the local flatness.}
We now consider a compact real $2$-codimensional submanifold $S$ of a
complex manifold $X$ of dimension $n\ge 3$.
\subsection{Induced foliation by CR orbits}
We use classical topological tools to obtain a description of the global structure
of the foliation by CR orbits.
\bp\Label{3.1}
Let $S\subset X$ be a compact connected real $2$-codimensional submanifold
such that the following holds:
\begin{enumerate}
\item[(i)] $S$ is nonminimal at every CR point;
\item[(ii)] every complex point of $S$ is flat and elliptic and there exists at least one
such point;
\item[(iii)] $S$ does not contain complex manifold of dimension
$(n-2)$. \end{enumerate}
Then $S$ is diffeomorphic to the unit sphere $\4S^{2n-2}\subset \C^{n-1}_z\times\R_x$
such that the complex points are the poles $\{x=\pm 1\}$
and the CR orbits in $S$ correspond to the $(2n-3)$-spheres given by $x=\const$. 
In particular, if  $S_{\rm ell}$ denotes the set
of all elliptic flat complex points of $S$, 
the open subset $S_0=S\setminus S_{\rm ell}$ 
carries a foliation $\mathcal F$ of class
$C^\infty$ with $1$-codimensional compact leaves.
\ep

\bpf
From conditions (i) and (ii), $S$ satisfies the hypotheses of Proposition 2.8. It follows from the proof of Proposition~\ref{elliptic-solution} 
that near an elliptic flat point,
all CR orbits are diffeomorphic to a real ellipsoid in $\R^{2n-2}$
(and hence to the sphere $\4S^{2n-3}$). 
Furthermore, the assumption (iii) guarantees that 
all CR orbits in $S$ must be of real dimension $2n-3$.
Hence, by removing small connected saturated neighbourhoods
of all complex points (since they are isolated, there are finitely many of them),
we obtain a compact manifold $S_0$ with boundary 
with the given foliation of codimension $1$ by its CR orbits whose first cohomology group with values in $\R$ is 0.

We now assume for a moment that this foliation is transversely oriented 
(see e.g.\ [CaC] for this and other terminology related to foliations).
Then we can apply Reeb-Thurston Stability Theorem 
(\cite{R,Th} and [CaC, Theorem~6.2.1]) to conclude
that $S_0$ is diffeomorphic to 
$\4S^{2n-3}\times [0,1]$ with CR orbits being of the form
$\4S^{2n-3}\times\{x\}$ for $x\in [0,1]$. 
Then the full manifold $S$ is  
diffeomorphic to the sphere $\4S^{2n-2}$ as required.

If $S_0$ with the foliation is not transversely oriented, we reach a
contradiction as   follows.
Consider its transversely oriented $2:1$ covering $\2S_0$.
Since the boundary leaves are simply connected and hence 
have trivial holonomies, $\2S_0$ has twice as  
many boundary leaves as $S_0$. As above, $\2S_0$ is diffeomorphic to  
$\4S^{2n-3}\times [0,1]$.
By continuity, each leaf $\4S^{2n-3}\times\{x\}$ of $\2S_0$ projects  
diffeomorphically onto a leaf
of $S_0$. Hence the base $[0,1]$ admits a $2:1$ covering over a  
$1$-dimensional manifold which
is impossible.
\epf
\subsection{Boundary value problem for embedded surfaces}
We now return to our central question:
{\it When does a compact
submanifold $S$ of $\C^n$
bound a Levi-flat hypersurface $M$?} 
From Proposition 3.1, we know that every CR
orbit of $S$ is a connected compact maximally complex CR submanifold of $\C^n$, $n\ge 3$,
and hence, in view of the classical result of
Harvey-Lawson [HaL], bounds a complex-analytic subvariety. 
Thus, in order to find $M$, at least as a real ``subvariety'',
foliated by complex subvarieties, 
a natural way to proceed is to build it as a family of
the solutions of the boundary problems for individual CR orbits. 
To do it, we reduce the problem to the corresponding problem 
in a real hyperplane of $\C^{n+1}$. 
The latter case is treated in the next section. 
It is inspired by the proof given in [Do1] for the $C^\omega$ case.

\section {On boundaries of families of holomorphic chains with $C^\infty$ parameters.}
As in [Do1] we follow the method of Harvey-Lawson [HaL] in [Ha, Section 3].
\subsection{Subvarieties with negligible singularities.}
Here we recall the terminology introduced in Definition~\ref{subvar}
with somewhat more details.
Let $X$ be a complex manifold endowed with a Hermitian metric. 
Let ${\6H}^d$ be the
$d$-dimensional Hausdorff measure on $X$. 
A closed subset $Y$ of $X$ is said to be a
{\it $d$-subvariety (with negligible singularities) of class $C^k$,
$k\in\N\cup\{\infty,\omega\}$},
if there exists a closed subset $\sigma$ of $Y$ such that ${\6
H}^d(\sigma)=0$ and $Y\setminus\sigma$ is a closed, {\sl oriented}
$d$-dimensional submanifold of class
$C^k$ of $X\setminus\sigma$ having 
locally (with respect to $X$) finite $\6 H^d$-measure.
%Similarly, a closed subset $Y\subset X$ is called
%a {\it $d$-subvariety with boundary (of class $C^k$)}
%if there exists a closed subset $\sigma$ of $Y$ such that 
%${\6H}^{d-1}(\sigma)=0$ and 
%$Y\setminus\sigma$ is a closed, {\sl oriented}
%$d$-dimensional submanifold with boundary (of class $C^k$)
%of $X\setminus\sigma$ having 
%locally (with respect to $X$) finite $\6 H^d$-measure.
The minimal set $\sigma$ as above is called the {\it
singular set} ${\rm Sing}\,Y$ of $Y$ 
(also called the {\it scar set} by Harvey-Lawson  
[Ha], [HaL1]) and
${\rm Reg}\,Y = Y\setminus\sigma$ its {\it regular part}. By
integration on ${\rm Reg}\, Y$,
we define a measurable, locally rectifiable current on $X$,
which is denoted by
$\lbrack Y\rbrack$
and said to
be the {\it integration current} of $Y$. 
The closed set $\sigma$ may be
increased without changing $\lbrack Y\rbrack$. The current 
$\lbrack Y\rbrack$ is said to be {\it closed} 
or to be a {\it cycle} if
$d\lbrack Y\rbrack=0$.
Obviously a change of the metric does not
change the class of $d$-subvarieties
as well as their regular and singular parts.

A $d$-subvariety $Y$ with negligible singularities (and the
current $\lbrack Y\rbrack$) is said to be CR of
CR dimension  $h$ 
if there exists a closed subset $\sigma'$ with 
$\sigma\subset\sigma'\subset Y$ 
such that ${\6
H}^d(\sigma')=0$
and $Y\setminus\sigma'$ is a CR submanifold of CR dimension $h$. 
In particular, $Y$ is said to be  {\it maximally complex}
or of {\it hypersurface type} if
$d=2h+1$ where $h$ is the CR dimension of $Y$; it is said to be 
{\it holomorphic} if it is a complex analytic variety.
A CR $d$-subvariety $Y$ of CR dimension $h$ 
is said to be {\it Levi-flat} 
if the regular part ${\rm Reg}\,Y$ is Levi-flat,
or equivalently, if ${\rm Reg}\,Y$ is foliated 
(in the usual sense) by complex $h$-dimensional submanifolds.
Furthermore, $Y$ is said to be 
{\it foliated by holomorphic $h$-varieties}, 
if through every point in $Y$ (including those of ${\rm Sing}\,Y$)
there is a complex $h$-dimensional subvariety in $Y$.

More generally, we call $d$-{\it chain} of $X$  every locally finite linear combination, with coefficients in $\Z$ of integration currents $[W_j]$ on $d$-subvarieties $W_j$ with negligible singularities. In particular, if the $W_j$ are complex subvarieties, the chain is said to be {\it holomorphic}.
%  whose closures ar

\subsection{Boundary problem in a real hyperplane of $\C^n$.}
\subsubsection{}  Here we extend to the $C^\infty$ case the $C^\omega$  
solution of the boundary  problem for a $C^\omega$ 
Levi-flat subvariety in a real  
hyperplane of $\C^n$ [Do1]. Let $n\geq 4$. 
We shall use the following
notation:  
$$z\sec=(z_2,\ldots,z_{n-1})\in\C^{n-2},
\quad\zeta'=(x_1,z\sec)\in\R\times\C^{n-2}.$$ 
Let $E=\R \times
\C^{n-1}=\{y_1=0\}\subset\C^n=\C\times\C^{n-1}$, and
$k\colon E\rightarrow
\R_{x_1}$, $(x_1;z\sec;z_n)\mapsto x_1$.
For $x^0_1\in \R_{x_1}$, set $E_{x_1^0}=k^{-1}(x_1^0)=\{x_1=x_1^0\}$.
\vskip 1mm
\subsubsection{} Let $N\subset E$ be a
compact, (oriented)  CR subvariety of $\C^n$ of
real dimension $2n-4$ and CR dimension
$n-3$, $(n\geq 4)$, of class $C^\infty$, with negligible
singularities (i.e. there exists a closed subset $\tau\subset N$ of
$(2n-4)$-dimensional Hausdorff measure $0$
such that $N\setminus \tau$ is a CR submanifold).
Let $\tau'$ be the set of all points $z\in N$
such that either $z\in\tau$ or $z\in N\setminus\tau$
and $N$ is not transversal to the complex hyperplane 
$k^{-1}(k(z))$ at $z$.
Assume that $N$,  as a
current of integration, is $d$-closed and satisfies:

\medskip
{\rm (H)} there exists a closed subset $L\subset\R_{x_1}$ 
with ${\6 H}^1(L)=0$ such that for every $x\in k(N)\setminus L$,
the fiber $k^{-1}(x)\cap N$ is connected 
and does not intersect $\tau'$.

\vskip 1mm

\subsubsection{}
We are going to prove the following:
\bt\Label{technical}
Let $N$ satisfy {\rm (H)} with $L$ chosen accordingly.
Then, there exists, in $E'= E\setminus k^{-1}(L)$, a unique $
C^\infty$ Levi-flat $(2n-3)$-subvariety $M$ with negligible singularities in $E'\setminus N$,
foliated by complex $(n-2)$-subvarieties, with the properties that
 $M$ simply (or trivially) extends  
to
$E'$ by a $(2n-3)$-current  (still denoted $M$) such that $dM=N$ in  
$E'$.
The leaves are the sections by the hyperplanes $E_{x_1^0}$, $x_1^0\in
k(N)\setminus L$, and are the solutions of the ``Harvey-Lawson problem'' for  
finding a holomorphic subvariety in $E_{x_1^0}\cong\C^{n-1}$ with prescribed boundary $N\cap E_{x_1^0}$.
\et

The proof will be given in the course of this section 4.\vskip 1mm

In what follows we increase $\tau$ such that $\tau=\tau'$ with $\tau'$ as above.
Recall that the {\it simple extension} $\widetilde M$ of the current $M$ in
$E'\setminus N$ to $E'$, if it
exists, satisfies, for every $\varphi\in {\6 D}(E')$, 
$<\widetilde M,\varphi>= <M,\varphi|_{E'\setminus N}>$.

In [Do1], [Do2], the statement [Do1, Th\'eor\`eme 6.9] is given for
$E'=E$, $n\geq 4$, and for $N$ being $C^\omega$-smooth, 
in two particular  
cases. Here we  shall give a proof for $C^\infty$ regularity, using  again  
[Ha], for any $N$  with negligible singularities outside  
$k^{-1}(L)$.

To be noted:

1) At the end of the proof, instead of the consideration of special
cases
as in [Do1], the technique of Harvey-Lawson [HaL] has been adapted.

2) In Theorem 4.1, the solution is a maximally complex $(2n-3)$-subvariety $M$ with negligible singularities and not a general chain as in [Do 1].

\subsection{Notations and reminders from [Do 1].}

\subsubsection{} Let 
$$
d\sec_E=\sum^n_{j=2}\frac{\partial}{\partial\overline
z_j}d\overline z_j, \quad
d'_E=\sum^n_{j=2}\frac{\partial}{\partial z_j}d z_j.
$$ 
In the following we choose an arbitrary 
orthogonal projection $\pi$ from $\C^n$ onto a complex $(n-1)$-plane.
Given $\pi$, after a linear coordinate change, we may assume that
$\pi:\C^n\rightarrow\C^{n-1}$, $z=(z_1,\ldots,z_n)\mapsto
(z_1,\ldots,z_{n-1})$, denotes the
projection as well as its restriction to $E$ or $E'$
(with $E'$ given in Theorem~\ref{technical}).  
Set ${\6 E}=\pi(E)$, ${\6 E}'=\pi(E')$, $\6N=\pi(N)$ and 
${\6 N}'=\pi(N\cap E')=\6N\cap \6E'$.
%Then ${\6 N}'=\pi(N\cap E')$ is a subvariety with negligible
%singularities, of codimension $1$ of ${\6 E}$.

Since $N$ is of CR dimension $n-3$, it follows:
$N=N^{3,1}+N^{2,2}+N^{1,3}$
in $\C^n$,
where $N^{r,s}$ denotes the component of type (bidegree) $(r,s)$.
Since $N$ is locally rectifiable, 
there exists a locally rectifiable current $P$ of
$E$ such that $j_*P=N$,
where $j\colon E\to \C^n$ is the inclusion. Then
\begin{equation}
P=P^{2,1}+P^{1,2}+dx_1\wedge(P^{2,0}+P^{1,1}+P^{0, 2}),
\end{equation}
where the types are relative to $(z_2,\ldots,z_n)$ and $dP=0$
(since $dN=0$).
In particular, taking the $(1,2)$ components of $dP=0$,
we obtain
\begin{equation}\Label{11}
d\sec_EP^{1,1}+d'_EP^{0,2}=\displaystyle
\frac{\partial P^{1,2}}
{\partial x_1}.
\end{equation}

\medskip
\subsubsection{} 
As in [Do 1] in $E$, we shall construct a 
{\it defining function} $R$ of our solution $M$ in $E'$ 
in the following way.
For $\zeta'\in {\6 E}'\setminus{\6 N}'$,
$z_n\in\C$, consider the Laurent series
\begin{equation}\Label{laurent}
\phi=C_0(\zeta') \ {\rm log} z_n+\sum_{m=1}^\infty m^{-1}C_m(\zeta')z_n^{-m},
\end{equation}
the coefficients $C_m$ being defined as follows.
Since $dP=0$, we have $d\sec_E P^{0,2}=0$.
Let $U^{0,1}$ be a solution
with compact support, of the equation $d\sec_EU^{0,1}=-P^{0,2}$
 (see \ [Do1, 5.4]). 
Set
\begin{equation}\Label{2}
C_m(\zeta')=K_E\sharp
   \pi_*\lbrack z_n^m(P^{1,1}+d'_EU^{0,1})\wedge dx_1\rbrack,
\end{equation}
where $K_E=\displaystyle\delta_0(z\sec)\otimes
H(x_1)\frac{\partial}{\partial x_1}$
with $H(x_1)$ being the Heaviside function in $x_1$,
and $K\sharp u$ is the 
{\it convolution-contraction} of a vector field 
$K=\sum K_j\frac{\partial}{\partial z_j}+
\sum K_{\bar j}\frac{\partial}{\partial \bar z_j}$
with distributional coefficients and a current 
$u=\sum u_{IJ}dz^I\wedge d\bar z^J$,
given by the formula
(cf. [HaL, p. 240] and [Ha, Section 3.6]):
$$K\sharp u = \sum 
(K_i * u_{IJ})
\,\imath_{\frac{\partial}{\partial z_i}}
(dz^I\wedge d\bar z^J)
+\sum 
(K_{\bar i} * u_{IJ})
\,\imath_{\frac{\partial}{\partial \bar z_i}}
(dz^I\wedge d\bar z^J),$$
where $\imath$ is the usual contraction.

Then each $C_m(\zeta')$ is a $(0,0)$-current and
by the construction, its restriction to
${\6 E}'\setminus{\6 N}'$ is represented by a $C^\infty$ function.
Furthermore, the series \eqref{laurent} converges  
for $(\zeta';z_n)\in({\6 E}'\setminus{\6 N}')\times(\C\setminus\overline{\Delta})$
and $R(\zeta';z_n)=\exp\phi(\zeta';z_n)$ is a $C^\infty$ function on 
$({\6 E}'\setminus{\6 N}')\times(\C\setminus\overline{\Delta})$, 
where $\Delta=\Delta (0,\rho)$ is a suitably large disc.

\subsection{Relation with slices} 
\bl\Label{slices}
Restricting to the hyperplane
$E_{x_1^0}\subset E$, for $x_1^0\in k(N)\setminus L\subset\R_{x_1^0}$, 
we obtain the hypotheses of the boundary problem of Harvey-Lawson in
$E_{x_1^0}\cong\C^{n-1}$ {\rm [Ha]}.
\el

Moreover, from condition (H), by Harvey-Lawson theorem [Ha, theorem 3.2], $N_{x_1^0}=k^{-1}(x_1^0)\cap N$ bounds an irreducible complex subvariety of $E_{x_1^0}\setminus N_{x_1^0}$ up to sign. 
\bpf 

Consider the slice $<C_m(\zeta'),k,x^0_1>$ 
of $C_m(\zeta')$ by the complex hyperplane $k^{-1}(x^0_1)$ 
(see e.g.\ [HaL]). 
From \eqref{2} we obtain
$$d\sec_EC_m(\zeta')=
d\sec_E\Big(K_E\sharp \pi_*\lbrack z_n^m(P^{1,1}+d'_EU^{0,1})\wedge
dx_1\rbrack\Big)=$$
$$d\sec_E\Big(\delta_0(z\sec)\otimes H(x_1)\frac{\partial}{\partial
x_1}\sharp \pi_*\lbrack
z_n^m(P^{1,1}+d'_EU^{0,1})\wedge dx_1\rbrack\Big)$$
$$=d\sec_E\Big(\big(\delta_0(z\sec)\otimes H(x_1)\big)* 
\big(\pi_*\lbrack z_n^m(P^{1,1}+d'_EU^{0,1})\rbrack\big)\Big)$$
$$=\big(\delta_0(z\sec)\otimes H(x_1)\big)*\big( \pi_*\lbrack
z_n^md\sec_E(P^{1,1}+d'_EU^{0,1})\rbrack\big)$$
As $d\sec_EU^{0,1}=-P^{0,2}$ by the construction in section 4.3.2, 
we obtain
$$
d\sec_E(P^{1,1}+d'_EU^{0,1})=d\sec_EP^{1,1}+d'_EP^{0,2}=\displaystyle
\frac{\partial P^{1,2}}
{\partial x_1}
$$ 
by \eqref{11} and
$$d\sec_EC_m(\zeta')=\Big(\delta_0(z\sec)\otimes \delta_0 (x_1)\Big)*
\Big(\pi_*\lbrack z_n^mP^{1,2}\rbrack\Big)
=\pi_*\lbrack z_n^mP^{1,2}\rbrack.$$ 
Denoting by $d\sec_0$ the $d\sec$
operator in $E_{x_1^0}$, we have
$$d\sec_0<C_m(\zeta'),k,x^0_1>=d\sec_E<C_m(\zeta'),k,x^0_1>=<d\sec_EC_m(\zeta'),
k,x^0_1>$$
$$=<\pi_*\lbrack z_n^mP^{1,2}\rbrack,k,x_1^0>,$$ \vskip 1mm

\noindent which is the equation satisfied by the coefficient 
$C_m(\zeta')$ of
Harvey-Lawson in $E_{x_1^0}$, for $N\cap  E_{x_1^0}$ [Ha, (3.1)]. 
By the uniqueness,
$<C_m(\zeta'),k,x^0_1>$ is the coefficient $C_m(\zeta')$ of Harvey-Lawson in $E_{x_1^0}$.
\epf

In view of Lemma~\ref{slices} and the properties of $C_m$'s
[HaL, \S6], $C_0(\zeta')$ is locally constant and integer-valued and
the function $R(\zeta';z_n)$ has the following properties:
\begin{enumerate}
\item[(a)] Let $V_0$ denote the union of the unbounded
connected components of 
${\6 E}'\setminus {\6 N}'$, 
then $R=1 \hskip 2mm on \hskip 2mm V_0\times
(\C\setminus\overline{\Delta})$;
\item[(b)] $R(\zeta';w)=\sum_{m=-\infty}^{C_0(\zeta')}A_m(\zeta')w^m,$
%\end{enumerate}
where the series converges uniformly on compact sets in 
$({\6E}\setminus{\6 N})\times
(\C\setminus\overline{\Delta})$;
%\begin{enumerate}
\item[(c)] every coefficient $A_m$ in (b) is a polynomial in a finite
number of $C_l$.
\end{enumerate}

Furthermore, for each fixed $x_1$, $R(x_1;\zeta'';z_n)$ extends
to $E_{x_1}\cap(({\6 E}'\setminus {\6 N}')\times\C)$ 
as a rational function in $z_n\in\C$ 
with holomorphic coefficients in 
$E_{x_1}\cap({\6 E}'\setminus {\6 N}')$.
The divisor of $R(x_1;\zeta'';z_n)$ gives the intersection
with $E_{x_1}\cap(({\6 E}'\setminus {\6 N}')\times\C)$
of the solution of Harvey-Lawson,
i.e.\ of the holomorphic chain in $E_{x_1}$ whose boundary (in the sense of currents) is either $E_{x_1}\cap N$ or $-E_{x_1}\cap N$.
Since $E_{x_1}\cap N$ is connected by our assumption (H) in 4.2.2, 
the extension of $R|_{x=x_1}$
is in fact either holomorphic or the inverse of a holomorphic function.
If $x_1\in k(N)\setminus L$, the slice $N\cap E_{x_1}$ is a nonempty 
compact maximally complex submanifold of $E_{x_1}$ of real codimension $3$ in view of (H).
Hence on every connected component of 
$k^{-1}(k(N))\cap(({\6 E}'\setminus {\6 N}')\times\C)$,
either $R(x_1;\zeta'';z_n)$ or its inverse is holomorphic
for fixed $x_1$. It then follows from the Cauchy integral formula
that the extension of $R$ (or of its inverse) 
is $C^\infty$ in all its variables.

\subsection{Construction of the maximally complex subvariety in
$\Gamma\setminus N$.}

We consider the projection $\pi :E\rightarrow {\6 E}$ and call
$\Gamma=\Gamma_\pi$ the set of
points $z\in E$ such that each point of $\pi^{-1}(\pi(z))\cap N$ is a
regular point of $\pi|_N$ (i.e.\ a smooth point of $N$ where
the differential of $\pi|_N$ is injective).

Denote by $(V_l)_l$ the family of the connected
components of ${\6 E}'\setminus{\6 N}'$.
Then for
$x_1^0\in k(N)\setminus L$  fixed, the function
$R(x_1^0;\zeta\sec;w)=R_{x_1^0}(\zeta\sec;w)$ is meromorphic in
$(\zeta\sec;w)$ on each $V_j\cap E_{x_1^0}$. From the Poincar\'e-Lelong
formula, 
$M_{x_1^0}=\displaystyle\frac{i}{\pi}\partial\overline\partial {\rm
log}|R_{x_1^0}(\zeta\sec;w)|$ is a holomorphic chain in the complement of $E_{x_1^0}\cap N$ in
$E_{x_1^0}\cap\Gamma$. In view of 4.4,  $\pm M_{x_1^0}$ is the
intersection with $\Gamma$ of the Harvey-Lawson solution
for $N\cap E_{x_1^0}$ in $E_{x_1^0}$. 
%Then, the $C^\infty$ coefficients that appear in
%$R(\xi_1,\zeta\sec;w)$ are not all flat for $\xi_1=\xi_1^0$.
Then locally $M_{x_1^0}=\displaystyle\sum m_l\lbrack Z_l\rbrack-\sum n_l\lbrack
P_l\rbrack$, where $Z_l$ (resp.\ $P_l$) are zero sets of the irreducible factors
$f_l$ (resp.\ $g_l$) of $R_{x_1}$ with multiplicities $m_l, n_l \geq 0$. 
In view of 4.4, for $\xi_1^0\in k(N)\setminus L$, the solution $M_{x_1^0}$
is the integration current on a complex subvariety $Z_{x_1^0}$
with multiplicity $m_{x_1^0}=\pm 1$ depending on the orientation 
of $N$. Define $M$ to be the union of $M_{x_1^0}$,
i.e.\ $M=\{(x_1,z''):z''\in M_{x_1}\}$.
Since $M$ is defined by $R$ (or $R^{-1}$), 
it is closed in $({\6 E}'\setminus{\6 N}')\times\C$.

\bl
 {\it Let $z^0$ be a point of $M$
belonging to a complex leaf $M_{x_1^0}$, where
$x_1^0=k(z^0)$. Assume that $z^0$ is smooth point of $M_{x_1^0}$.
Then there exists a
neighbourhood $U_{z^0}$ of $z^0$ in $\C^n\setminus N$ such that} $M\cap U_{z^0}$ {\it is smooth. Moreover, the multiplicity $m=m_{x_1^0}$ is locally constant.}
\el
\bpf
Consider the defining function $R(x_1;\zeta\sec;z_n)$ of $M$
in the neighbourhood of $z^0=(x_1^0;{\zeta\sec}^0;z_n^0)$. The
defining function of $M_{x_1^0}$ is $R(x_1^0;\cdot;\cdot)$. 
It is rational in $z_n$ with coefficients being smooth functions in $(x_1,\zeta\sec)$ in view of 4.4. From the regularity in $z^0$, in the above notations, we have $M_{x_1^0}= m_{x_1^0}\lbrack Z_{x_1^0}\rbrack$, $m_{x_1^0}=\pm 1$, and $R(x_1;\zeta\sec;z_n)=(f(x_1;\zeta\sec;z_n))^{\pm1}$, where $f(x_1;\zeta\sec;z_n)$ is $C^\infty$, holomorphic in $(\zeta\sec;z_n)$ and irreducible at $z^0$. Hence the gradient of $f(x_1;\cdot;\cdot)$ does not vanish at $({\zeta\sec}^0;z_n^0)$.
Then locally near $z_0$, $M$ is given by $f=0$ and is therefore a $C^\infty$ submanifold
by the implicit function theorem. 
%$(\xi_1^0,\zeta\sec;z_n)$
Let $Z=f^{-1}(0)$ locally; then $Z_{\xi_1^0}=Z\vert_{E_{\xi_1^0}}$, and there exists a neighbourhood $U_{z_0}$ of $z_0$ on which $Z$ is smooth. The multiplicity $m$ is constant in view of 4.4.
\epf
\noindent {\bf Corollary}.- {\it The integer $m$ is constant on each
connected component of $(k(N)\setminus L)\times\C^{n-1}$.}
$\square$\vskip 1mm

In view of Lemma 4.3, the set of the singular points of $M$ is the union of the singular points of 
the complex analytic sets $M_{x_1}$ for $x_1$ in a connected component of $k(N)\setminus L$, since
the complement of this set is smooth. The set of singular points of each $M_{x_1}$ is of real dimension at
most $2n-6$, $x_1\in\R$, consequently the Hausdorff dimension of the singular set of $M$ is at most
$2n-5$.

For $x_1^0\in L$, the situation is unknown at the moment.

The defining function $R$ and the variety $M$ are defined in 
${\big((\6E'\setminus \6N')\times
\C\big)}\cap E'$. To extend $M$ into a maximally complex $(2n-3)$-variety
in $(\Gamma\setminus N)\cap E'$, we have to extend $M$ to any point $z_0\notin N$ with 
$\pi(z_0)=\zeta'\in\6N'$. That is done by solving the boundary problem in the
neighborhood of $z_0$ (for details, see  [Do1, Proposition 6.6.2], [Ha], Lemma 3.22).

The variety $M$ of $\Gamma\setminus N$ has a finite volume in the neighborhood of every point 
$z^0\in N\cap\Gamma$ because, for $x_1$ in the neighborhood of $k(z^0)$ in $k(N)\setminus L$, $M_{x_1}$
is of finite $ (2n-4)$-volume, so $M$ is of finite $(2n-3)$-volume. Then
$M$ has a simple extension, still denoted $M$, to a
rectifiable current of dimension $(2n-3)$ in $\Gamma\cap E'$. 

Moreover, on every connected component $E'_0$ of $E'$, $dM=\displaystyle\sum n_j\lbrack N_j\rbrack$, 
where the $N_j$ are the connected components of $N\cap\Gamma\cap E'_0$, from classical properties of
locally rectifiable currents [Do1, Proposition 6.7.2].

\subsection{End of the proof.}
We now give a sketch of the end of the proof of Theorem 4.1.
Let $\Gamma'=\Gamma\cap E'$.
By a generic choice of the projection $\pi$,
we may assume that ${\6H}^{2n-4}(N\setminus\Gamma)=0$.

\bp\Label{bdry}
There exists a closed subset $A$ of $N$ 
(containing $N\setminus \Gamma'$) such
that ${\6 H}^{2n-4}(A)=0$
and such that each point 
$z^0=(x_1^0,z''{}^0, z_n^0)$
of $(N\cap\Gamma')\setminus A$
has a small open neighbourhood $B$ with the following properties:
\begin{enumerate}
\item[a)] $B=I\times B\sec\times\Delta$, where $I$ is an interval of
$\R_{x_1}$ centered at $x_1^0$,
$B\sec$ a ball of $\C^{n-2}_{z\sec}$ centered at ${z\sec}^0$, and $\Delta$
a disc of
$\C_{z_n}$ centered at $z_n^0$;\\
\item[b)] $W={\rm supp} \ M\cap (B\setminus N)$ is a submanifold, 
and
either
\begin{enumerate}
\item[(i)] $W$ is connected and the pair $(\pm W, N\cap B)$ is a
$C^\infty$
submanifold with boundary, or\\
\item[(ii)] $\overline W =W\cup (N\cap B)$ is a connected, maximally
complex
submanifold of $B$ containing $N\cap B$ as a real $C^\infty$
hypersurface, dividing
$\overline W$ into two components $W_i$ and $W_j$	over  $(I\times
B\sec)\cap V_i$ and
$(I\times B\sec)\cap V_j$ in ${\6 E}$, respectively.
\end{enumerate}
\end{enumerate}
In either case, if the labelling is chosen so that $\pi(N\cap B)$ is
the oriented boundary of $V_i\cap (I\times B\sec)$, then $M\mid_B$ is of the form
$(m+1)[W_i]+m[W_j]$ so that $dM=N$ on $B$. ({\it In case $(i)$, $m=-1$ or $m=0$}).
Moreover, $A$ can be chosen such that, for every 
point $z\in N'\cap E_x$ where $\pi|_N$ is regular,
one has $\6H^{2n-3}(A\cap E_x\cap U_z)=0$
for some neighborhood $U_z$ of $z$.
\ep

%Proposition 4.4 contains a lemma for what follows and a result on  
%singularities on the boundary in addition to Theorem 4.1.
\bpf
The proof proceeds as the one of Lemma 3.24 of $\lbrack$Ha$\rbrack$ for
fixed $x_1\notin L$ (see also [HaL]). Suppose $B\subset V_j\times \C$ as in
 Proposition~4.4. Consider the
function $R_j(x_1;z\sec;w)$ defined on $V_j\times \C$. In $B$,  $M$ is defined by the vanishing of the function $p_j(x_1,z\sec,w)$,
product of the distinct factors of the numerator and the denominator of  
the
rational function $R_j$ in $w$. The function $p_j(x_1,z\sec,w)$ can be
considered as a unitary polynomial in $w$, up to a non vanishing factor. Since
the coefficients of $p_j$ are $C^\infty$,
$B$ being contained in a connected component of  
$((k(N)\setminus L)\times\C^{n-1})\cap\Gamma$, we get the conclusion.
\epf

\subsubsection{Varying the projection $\pi$.} (cf. [Do1], 5.1.1, 6.1.2.)

\bl
Let $\pi, \pi'\in\cp {n-1}$ be two projections,
$\Gamma'_{\pi}, \Gamma'_{\pi'}$ be the corresponding sets $\Gamma'$.  
Then the currents $M_{\pi}$ in $\Gamma'_{\pi}$, $M_{\pi'}$ in $\Gamma'_{\pi'}$agree on
$\Gamma'_{\pi}\cap\Gamma'_{\pi'}$, and their union is the current $M$
in $E'\setminus N$ with $dM=N$.
\el

\begin{proof}(From the proof of [Ha, Lemma 3.25].)
It suffices
to prove the first part of  the conclusion for $\pi$ and $\pi'$ close
enough. It is known and easy to show (see [Do1]) that $M=j_*S$ where $S$
is a locally rectifiable current of $E'$ and $j:E'\rightarrow \C^n$ is the inclusion and where the
coefficient $C_m(\zeta')=\pi_*(z_n^mS_\pi)$ on $\Gamma'_\pi$; the projection $\pi$ and hence the 
coordinate $z_n$
being chosen, let $C'_m(\zeta')=\pi_*(z_n^mS_{\pi'})$. Then $\displaystyle\big(d\sec
_E+\frac{\partial}{\partial x_1}dx_1 \big) C_m=\pi_*(z_n^m P^{1,2})=\displaystyle\big(d\sec
_E+\frac{\partial}{\partial x_1}dx_1 \big) C'_m$. The difference $c_m=C_m-C'_m$  satisfies the equation
$\displaystyle\big ( d\sec _E+\frac{\partial}{\partial x_1}dx_1 \big )  
c_m=0$. Let $B$ be a closed ball of $E$ containing $M$, then it can be shown that $C_m$ and $C'_m$ are
defined on an open set $U'$ of ${\6E}'$ containing points ouside $\pi(B)$ and therefore $C_m$ and $C'_m$
vanish on $U'\setminus\pi(B)$; since $c_m$ is holomorphic in $z''$ and constant in $x_1$, one has
$c_m=0$ on $U'$.

\epf

This finishes the proof of Theorem~\ref{technical}.

\section{On the existence of a Levi-flat $(2n-1)$-subvariety with prescribed boundary.}
We now return to the initial problem of finding a real Levi-flat
hypersurface in $\C^n$ with
prescribed boundary. We translate this problem into a boundary problem for subvarieties of a
hyperplane $E$ of $\C^{n+1}$ with negligible singularities, 
foliated by holomorphic chains and then apply Theorem 4.1.
We mention that Delannay [De] gives a solution of the problem under  
certain additional assumptions.
\bpf[Proof of Theorems~\ref{main} and \ref{precise}]
We first show that there exists a global smooth function $\nu\colon S\to [0,1]$ without critical points away from the complex points such that the complex points are $\nu^{-1}(j)$, $j=0,1$, and the level sets of $\nu$ are the CR orbits of $S$ away from the complex points.
By Proposition~\ref{elliptic-solution}, such $\nu$ can be constructed near 
every complex point. Furthermore, in view of Proposition~\ref{3.1},
such $\nu$ can be obtained globally on $S$ away from its complex points.
Putting everything together and using a partition of unity,
we obtain a function $\nu\colon S\to [0,1]$ with desired properties.

The submanifold $S$ being, locally, a boundary of Levi-flat hypersurface, is orientable. We now set $\2S=N={\rm gr}\,\nu = \{(\nu(z),z) : z\in S\}$. The corresponding map $\lambda\colon S\rightarrow\2S$, $z\mapsto \nu ((z),z)$ is diffeomorphic; moreover $\lambda\vert_{S\setminus S_{\rm ell}}$ is a CR map. Choose an orientation on $S$. Then $N$ is an (oriented) CR subvariety with the two-point set of singularities $\tau=\lambda(S_{\rm ell})$.

At every point of $S\setminus S_{\rm ell}$, $d_{x_1} \nu\not = 0$, then condition (H) of section 4.2.2 is satisfied with $L=k(S_{\rm ell})$.
Then all the assumptions of Theorem~\ref{technical}
are satisfied. We conclude that  $N$ is the boundary 
of the Levi-flat $(2n-1)$-subvariety $\2M$ in $[0,1]\times\C^n\subset\R\times\C^n=E$. The desired conclusions now follow, at least away from the complex points, from Theorems~\ref{technical} and \ref{bdry}, where we take the intersection of the sets $A$ in Theorem~\ref{bdry} for all projections $\pi$. The description near complex points is given in Remark~\ref{filling}.
\epf

\section{On the uniqueness of a Levi-flat $(2n-1)$-subvariety
with prescribed boundary.}

Here we shall prove Theorem~\ref{unique}.
We begin by establishing an elementary auxiliary lemma.
Recall that a domain with $C^2$-smooth boundary is 
{\em strongly convex} if its boundary is locally at each point given
by vanishing of a $C^2$-smooth function whose 
Hessian is positive definite.

\bl\Label{include}
Let $S\subset\R^m$ be a compact, smooth submanifold.
Then there exists a bounded strongly convex domain in $\R^m$
with $C^2$-smooth boundary whose closure contains $S$ and whose boundary
contains a nonempty open subset of $S$.
\el

\bpf
Let $B$ be the minimal open ball in $\R^m$ with center $0$ whose closure contains $S$.
Then there exists a point $p\in \partial B\cap S$. Without loss of generality,
the ball is unit and $p=(0,\ldots,0,-1)$. Furthermore, we may assume that 
$T_pS=\R^d\times\{0\}$ with $d=\dim S$.
Near $p$, the ball $B$ is given by $x_m>-1+\rho_0$, 
where 
$$\rho_0(x_1,\ldots,x_{m-1}):=1-\sqrt{1-x_1^2-\ldots-x_{m-1}^2}.$$
The submanifold $S$ near $p$ is given by 
$x_j=\phi_j=\phi_j(x_1,\ldots,x_d)$ for $j=d+1,\ldots,m$,
for a suitable smooth functions $\phi_j$.
Since $S\subset \1B$, we have for all $(x_1,\ldots,x_d)$ near $0$,
$$\phi_m(x_1,\ldots,x_d)\ge -1 + 
\rho_0\big(x_1,\ldots,x_d,\phi_{d+1}(x_1,\ldots,x_d),\ldots,
\phi_{m-1}(x_1,\ldots,x_d)\big),$$
 and therefore $\rho(x_1,\ldots,x_{m-1})\ge\rho_0(x_1,\ldots,x_{m-1})$, where
$$\rho(x_1,\ldots,x_{m-1}):= \rho_0(x_1,\ldots,x_{m-1}) +
\phi_m(x_1,\ldots,x_d)+1 - 
\rho_0\big(x_1,\ldots,x_d,\phi_{d+1}
(x_1,\ldots,x_d),\ldots,\phi_{m-1}(x_1,\ldots,x_d)\big).$$
Furthermore, we clearly have $x_m=\rho(x_1,\ldots,x_{m-1})$ 
whenever $x\in S$ and the function $\rho$ is strongly convex.
It remains to ``glue'' the functions $\rho_0$ and $\rho$,
i.e.\ to construct a smooth strongly convex function
$\2\rho$ which coincides with $\rho$ in a neighborhood $U_1$ of $0$
and with $\rho_0$ outside another neighborhood $U_2\supset U_1$.
Thus the proof is completed by applying in radial directions
the following elementary one-dimensional lemma whose proof is omitted.
\epf

\bl
Let $\rho_0(x)\le\rho(x)$ be two strongly convex 
$C^2$-smooth functions in $x\in\R$.
Then there is a third strongly convex 
$C^2$-smooth function $\2\rho(x)$ with 
$\rho_0(x)\le\2\rho(x)\le\rho(x)$ and two neighborhoods
$U_1\subset U_2$ of $0$,
such that $\2\rho=\rho$ in $U_1$ and $\2\rho=\rho_0$ outside $U_2$.
The neighborhoods $U_1$ and $U_2$ can be chosen arbitrarily small.
Moreover, if $\rho$ and $\rho_0$ depend $C^2$-smoothly on an additional
parameter $y$ that belongs to a compact manifold $K$,
then $U_1$ and $U_2$ can be chosen arbitrarily small and
uniform for all $y\in K$ and $\2\rho$ can be chosen
$C^2$-smooth in $(x,y)$.
\el

Recall that by a leaf in a Levi-flat 
$(2n-1)$-subvariety $M$ of $\C^n$
we mean a maximal connected immersed complex hypersurface in 
${\rm Reg}\, M$, i.e.\ a connected complex $(n-1)$-dimensional 
manifold $\2L$ with an injective holomorphic immersion 
$j\colon \2L\to M$ such that $L=j(\2L)$ cannot be enlarged.
By the {\em intrinsic topology} on $L$ we mean the one induced by $j$.

\bpf[Proof of Theorem~\ref{unique}]
Let $\Omega\subset \C^n$ be a bounded strongly convex domain
satisfying the conclusion of Lemma~\ref{include}.
In view of (ii) (where we always refer to Theorem~\ref{unique}), 
$M$ has to be contained in $\1\Omega$
by the maximum principle. By the same argument,
also the image $\pi(\2M)$ has to be contained in $\1\Omega$.

Let $A\subset S$ be as in Theorem~\ref{unique} (iii)
and $\2A\subset\2S$ be as in Theorem~\ref{precise} (iii).
Since $S\cap \d\Omega$ 
has nonempty interior,
we can find a point $p\in (S\cap \d\Omega)\setminus A$.
In addition, we may assume that $p$ is a CR point of $S$.
Then near $p$, $M$ is a Levi-flat hypersurface with boundary $S$ by (iii).
Since $p$ is CR, the complex leaves $L$ of a sufficiently
small neighborhood of $p$ in $M$
are closed complex submanifolds with CR orbits of 
$S$ as their boundaries.
On the other hand, the same holds for the images $\pi(\2L)$
of the leaves $\2L$ of $\2M$ away from $\2A$.
Recall from Theorem~\ref{precise} (i) that a leaf $\2L$ of $\2M$
is the regular locus of a complex-analytic $(n-1)$-dimensional
subvariety $(\2M\setminus\2S)\cap (\{x\}\times\C^n)$ for
$x\in (0,1)$. We shall write $\2L_x$ to indicate the dependence on $x$.
Without loss of generality, 
$\2p:=(\pi|_{\2S})^{-1}(p)\notin \2A$.
Since both $L$ and $\pi(\2L)$ are in $\1\Omega$
and $\Omega$ is strongly convex,
it follows that, whenever $L$ and $\pi(\2L)$
are bounded by the same CR orbit, they must coincide
in a neighborhood of $p$
(we have used the well-know fact that a holomorphic function
is uniquely determined by its boundary value on a real hypersurface).
Consequently, a neighborhood of $p$ in $M$ coincides
with the image of a neighborhood of $\2p$ in $\2M$.

We next consider any global complex leaf $\2L=\2L_{x_0}$ of $\2M$
with $x_0\in (0,1)$ such that $M$ contains 
a nonempty open subset of $\pi(\2L)$. 
Assume first that $\6H^{2n-3}({\rm Sing}\, M\cap \pi(\2L))=0$.
Then $\pi(\2L)\setminus {\rm Sing}\, M$ is connected.
Since ${\rm Reg}\, M$ is a Levi-flat hypersurface
and $\pi(\2L)$ is a (locally closed) smooth 
complex hypersurface in $\C^n$,
it is easy to see that the set of all points of
$\pi(\2L)\cap M$ that are interior with respect to $\pi(\2L)$
is both open and closed in $\pi(\2L)\setminus {\rm Sing}\, M$.
This shows that $\pi(\2L)\setminus {\rm Sing}\, M\subset M$
and therefore $\pi(\2L)\subset M$.

On the other hand, suppose that 
$\6H^{2n-3}({\rm Sing}\, M\cap \pi(\2L))>0$.
We claim that this can only happen for at most
countable set $\6S$ of leaves $\2L$.
Indeed, this follows immediately from (i)
and the fact that $\dim(\pi(\2L_1)\cap \pi(\2L_2))\le 2n-4$
whenever $\2L_1$ and $\2L_2$ are two distinct leaves.

We are now ready to prove that $\pi(\2L_{x_0})\subset M$ for any $x_0$.
Indeed, it follows from the above, that the set $\6A$ of such $x_0$
is nonempty. Since $M$ is closed, the set $\6A$ is also closed
in the interval $(0,1)$. We have to show that it is open.
Fix any $x_0\in \6A$ and consider $\2L=\2L_{x_0}$.
Since $\dim \2L=2n-2$, we have $\pi(\2L)\not\subset {\rm Sing\, M}$
in view of (i).
Since each connected component of 
$L:=\pi(\2L)\setminus {\rm Sing\, M}\subset {\rm Reg}\, M$
is a subset of a leaf of $M$, 
we conclude from (iii) that $\1{\pi(\2L)}\cap S$
is not contained in $A$.
Then we can find a point 
$q\in (\1{\pi(\2L)}\cap S)\setminus A$ such that 
$\2q:=(\pi|_{\2S})^{-1}(q)\notin \2A$ (where we use 
the conclusion (iii) of Theorem~\ref{precise} that is assumed
to be satisfied).
Then both $M$ and the image of a neighborhood of $\2q$ in $\2M$
are Levi-flat hypersurfaces with the same boundary $S$.
Moreover, they both contain an open subset of $\pi(\2L)$.
Arguing as above, we see that these two Levi-flat hypersurfaces
must coincide near $q$. In particular, for any leaf $\2L_{x_1}$
with $x_1$ sufficiently close $x_0$, 
$M$ contains a nonempty open subset of $\pi(\2L_{x_1})$.
Above we have seen that either $\pi(\2L_{x_1})\subset M$
or $\2L_{x_1}$ is one of the countably many leaves
with $\6H^{2n-3}({\rm Sing}\, M\cap \pi(\2L_{x_1}))>0$.
In the latter case, $x_1$ is a limit
of a sequence of points $x'_1$ with 
$\6H^{2n-3}({\rm Sing}\, M\cap \pi(\2L_{x'_1}))=0$.
Then $\pi(\2L_{x'_1})\subset M$ for each such $x'_1$
and therefore $\pi(\2L_{x_1})\subset M$ because $M$ is closed.
This shows that $\6A$ is indeed open and thus completes the proof
that $\pi(\2L_{x_0})\subset M$ for any $x_0$.

Thus $\pi(\2M)\subset M$ by taking the closure.
We finally claim that $\pi(\2M)=M$.
It suffices to show that ${\rm Reg}\, M\subset \pi(\2M)$.
Assume by contradiction that there exists 
$z\in {\rm Reg}\, M \setminus \pi(\2M)$ 
and consider the leaf $L\subset {\rm Reg}\,M$ through $z$.
In view of (iii), there exists a point 
$q\in (\1L\cap S)\setminus A$, which is a CR point $S$.
Since near $q$, $M$ is a Levi-flat hypersurface with boundary $S$,
it follows that $L$ contains a complex hypersurface 
with boundary being one of the CR orbits of $S$.
On the other hand, there exists a (unique) leaf $\2L$ 
of $\2M$ such that
$\pi(\2L)$ contains $q$ in its closure.
In fact, $\2L=\2L_x$ with $x$ being the projection to $\R$ of
$\2q:=(\pi|_{\2S})^{-1}(q)\in \R\times\C^n$.
Moreover we have shown that $\pi(\2L)\subset M$.
Then $\pi(\2L)$ must contain an open subset of $L$
(with respect to the intrinsic topology)
and therefore $\pi(\2L)$ contains the whole $L$.
In particular, $\pi(\2L)$ contains the point $z$,
which leads to a contradiction. The proof is complete.
\epf

%\newpage

\end{document}